\documentclass [11pt]{amsart}
\usepackage {amsmath, amssymb, amscd, graphicx, color, url}
\newtheorem {theorem}{Theorem}

\newtheorem {lemma}[theorem]{Lemma}

\newtheorem {corollary}[theorem]{Corollary}
\newtheorem {conjecture}[theorem]{Conjecture}

\def\zz {{\mathbb{Z}}}
\def\rr {{\mathbb{R}}}

\def\ff {{\mathbb{F}}}

\def\hh {{\mathcal{H}}}
\def\tt {{\mathbb{T}}}
\def\xx {{\mathbf{x}}}
\def\yy {{\mathbf{y}}}
\def\del {{\partial}}
\def\kh {{\widetilde{Kh}}}
\def\hfk {{\widehat{HFK}}}
\def\hfd {{\widehat{HF}(-\Sigma(L))}}

\def\rk {{\operatorname{rk}_{\ff}}}
\def\wsigma {{\widehat {\Sigma}}}
\def\sym {{\mathrm{Sym}}}
\def\ta {{\tt_{\alpha}}}
\def\tb {{\tt_{\beta}}}
\def\tg {{\tt_{\gamma}}}
\def\td {{\tt_{\delta}}}
\def\tbp {{\tt_{\beta'}}}
\def\fin {{\hfill \square}}

\def\det {{\operatorname{det}}}

\def\To {{\ \rightarrow \ }}
\begin{document}

\title
[An unoriented skein exact triangle]{An 
unoriented skein exact triangle for knot Floer homology}

\author [Ciprian Manolescu]{Ciprian Manolescu}
\thanks {The author was supported by a Clay Research Fellowship.}
\address {Department of Mathematics, Columbia University\\ New York, NY 10027}
\email {cm@math.columbia.edu}

\begin {abstract} 
Given a crossing in a planar diagram of a link in the three-sphere, we show 
that the knot Floer homologies of the link and its two resolutions at that 
crossing are related by an exact triangle. As a consequence, we deduce that 
for any quasi-alternating link, the total rank of its knot Floer homology 
is equal to the determinant of the link. 
\end {abstract}

\maketitle

\section {Introduction}

Given a link $L \subset S^3,$ there are various interesting {\it knot 
homology theories} associated to $L.$ One of them is Khovanov's reduced 
theory $\kh(L)$ from \cite{K},\cite{K2}, a bigraded vector space over $\ff 
= \zz/2\zz$ whose Euler characteristic is the Jones polynomial. Other 
theories can be obtained from symplectic geometry: for example, one can 
consider $-\Sigma(L),$ the double cover of $S^3$ branched along $L$ (with 
the orientation reversed), and apply to it the Heegaard Floer homology 
$\widehat{HF}$ functor of Ozsv\'ath and Szab\'o \cite{OS1}. In a similar 
vein, a very useful theory is the knot Floer homology $\hfk(L)$ of 
Ozsv\'ath-Szab\'o and Rasmussen \cite{OS2}, \cite{R1}. In its simplest 
form, $\hfk(L)$ is a bigraded vector space whose Euler characteristic is 
the Alexander polynomial. Knot Floer homology is known to detect the genus 
of a knot \cite{OS3.5}, as well as whether a knot is fibered \cite{N}. 
There exists a refinement of $\hfk$ called link Floer homology \cite{OS5}, 
which detects the Thurston norm of the link complement \cite{OS6}. For the 
purposes of this paper, we will consider all the theories with 
coefficients in $\ff = \zz/2\zz,$ and usually ignore all their gradings.

The three theories mentioned above are rather different in origin. 
Khovanov's theory was developed out of a study of representation theory 
and categorification, and was defined in a purely combinatorial fashion to 
start with. On the other hand, the Floer homologies $\hfd$ and $\hfk(L)$ 
were originally constructed using pseudolomorphic disks, and a 
combinatorial description of them has only recently been found 
(\cite{MOS}, \cite{SW}). The ways in which these three theories are 
related to one another, however, are still not completely understood.

\begin {figure}
\begin {center}
\input {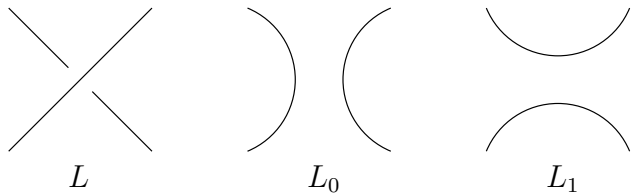}
\caption {The links in the unoriented skein relation.}
\label {fig:skein}
\end {center}
\end {figure}

One connection between $\kh(L)$ and $\hfd$ was pointed out by Ozsv\'ath 
and 
Szab\'o in \cite{OS4}. They observed that both of these theories satisfy  
unoriented skein exact triangles. More precisely, consider three links 
$L, L_0$ and $L_1$ that admit planar diagrams differing from each other only 
at one crossing, where they look as in Figure~\ref{fig:skein}. Then there is 
an exact triangle
$$ \hh(L) \To \hh(L_0) \To \hh(L_1) \To \hh(L), $$ 
where the symbol $\hh$ could stand for either $\kh$ or $\widehat{HF} 
(-\Sigma(\cdot)).$ One can iterate this exact triangle by 
applying it to all crossings of $L.$ Via some homological algebra, the 
result of this iteration process is a spectral sequence whose $E^2$ term 
is Khovanov homology, and which converges to  $\hfd.$ This implies an 
inequality of ranks: 
\begin {equation}
\label {ineq}
\rk \kh(L) \geq \rk \hfd. 
\end {equation}
In \cite{OS4}, Ozsv\'ath and Szab\'o defined a class of links, called 
quasi-alternating links, for which (\ref{ineq}) becomes equality and, 
furthermore, the ranks of the two theories are equal to $\det(L),$ the 
determinant of the link. In particular, all alternating links are 
quasi-alternating.

The main result of this paper is that knot Floer homology also has an
unoriented skein exact triangle. Note that it was previously known 
to satisfy an oriented skein exact triangle, cf. \cite{OS2}.

\begin {theorem}
\label {main}
Let $L$ be a link in $S^3.$ Given a planar diagram of $L$, let $L_0$ and  
$L_1$ be the two resolutions of $L$ at a crossing in that diagram, as in  
Figure~\ref{fig:skein}. Denote by $l,l_0, l_1$ the number of 
components of the links $L, L_0,$ and $L_1,$ respectively, and set $m = 
\max \{ l, l_0, l_1\}.$ Then, there is an exact triangle
$$
 \hfk(L) \otimes V^{m-l} \To \hfk(L_0)  \otimes V^{m-l_0} \To \hfk(L_1) 
 \otimes V^{m-l_1} \To \hfk(L) \otimes V^{m-l}, $$
where $V$ denotes a two-dimensional vector space over $\ff.$
\end {theorem}

Ozsv\'ath and Szab\'o \cite{OS3}, \cite{OS5} proved that $\rk \hfk(L) = 
2^{l-1} \det (L)$ when $L$ is an alternating link with $l$ 
components. A simple consequence of Theorem~\ref{main} is a generalization 
of their result:

\begin {corollary} \label {cor} 
If $L$ is a quasi-alternating link with 
$l$ components, then $\rk \hfk(L)  = 2^{l-1} \cdot \det(L).$ 
\end {corollary}

Rasmussen \cite[Section 5]{R2} observed that $\kh$ and $\hfk$ have equal 
ranks for many classes of knots (including most small knots), and 
asked for an explanation. Corollary~\ref{cor} can be viewed as a partial 
answer to his question. Indeed, as we show in Section~\ref{sec:qa}, many 
small knots are quasi-alternating. 

As suggested in \cite{R2}, a more convincing explanation for the 
similarities between $\kh$ and $\hfk$ (for knots) would be a spectral 
sequence whose $E^2$ term is $\kh$ and which converges to $\hfk.$ (For 
links of $l$ components, $\kh$ should be replaced with $\kh \otimes 
V^{l-1}$.) This would imply an inequality of ranks similar to 
(\ref{ineq}), namely $2^{l-1} \cdot \rk \kh \geq \rk \hfk.$ In turn, using 
the fact that $\hfk$ detects the unknot \cite{OS3.5}, this would provide a 
positive answer to the following well-known conjecture:

\begin {conjecture}
\label {unknot}
If $K$ is a knot with $\rk \kh(K) = 1,$ then $K$ is the unknot.
\end {conjecture}

Unfortunately, Theorem~\ref{main} does not directly imply the inequality 
of ranks and hence Conjecture~\ref{unknot}. One can iterate the unoriented 
skein triangle for $\hfk$ and obtain a spectral sequence whose 
$E^{\infty}$ term is $\hfk \otimes V^n,$ for some large $n.$ However, the 
presence of the $V$ factors makes it unclear whether the $E^2$ term of the 
spectral sequence is related to Khovanov homology. We leave this as an 
open problem.

We end this introduction with some remarks. First, the exact triangle from 
Theorem~\ref{main} is different from the other exact triangles in Floer 
homology, in that the maps do not even respect the homological gradings 
modulo $2.$ Also, the proof only works for the hat version of $HFK,$ and 
it is not clear whether a similar triangle holds for the other versions. 
Finally, in this paper we only use the traditional definition of knot 
Floer homology, based on counting pseudo-holomorphic disks. It would be 
interesting to recover the same result, and perhaps do further 
computations, using the combinatorial definition of \cite{MOS} instead.

\medskip \noindent {\bf Acknowledgments.} I owe a great debt of gratitude 
to Peter Ozsv\'ath; his suggestions and encouragement have been essential 
in completing this paper. I am also grateful to Nathan Dunfield and Jacob 
Rasmussen for several very helpful discussions, and John Baldwin for 
pointing out several minor errors in a previous version.

\section {Special Heegaard diagrams}
\label {sec:special}

For the original definitions of knot Floer homology, we refer the reader 
to \cite{OS2}, \cite{R1}, \cite{OS5}, \cite{MOS}. Here we will use 
a special class of Heegaard diagrams, which are particular cases of  
the multiply-pointed diagrams defined in \cite[Section2]{MOS}. They are a 
variant of the Heegaard diagrams asssociated to knot projections in 
\cite[Section 2]{OS3}. 

Let $L \subset S^3$ be a link with $l$ components and $D$ a planar, 
connected projection of $L.$ If $D$ has $c$ crossings, then it splits the 
plane into $c+2$ regions. Let $A_0$ denote the unbounded region in $\rr^2 
- D,$ and $A_1$ a region adjacent to $A,$ separated from $A$ by an edge 
$e.$ Denote the other regions by $A_2, \dots, A_{c+1}.$ Let $s_1, \dots, 
s_k \ (k \geq l-1)$ be a collection of (not necessarily distinct) edges of 
$D,$ chosen in such a way that for every component of the link, its 
projection contains at least one of the edges $s_i$ or $e.$

\begin {figure}
\begin {center}
\begin{picture}(0,0)%
\includegraphics{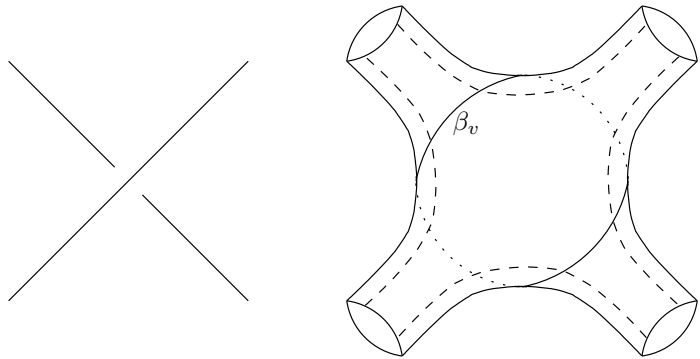}%
\end{picture}%
\setlength{\unitlength}{3947sp}%
\begingroup\makeatletter\ifx\SetFigFont\undefined%
\gdef\SetFigFont#1#2#3#4#5{%
  \reset@font\fontsize{#1}{#2pt}%
  \fontfamily{#3}\fontseries{#4}\fontshape{#5}%
  \selectfont}%
\fi\endgroup%
\begin{picture}(4349,2224)(-2261,-1973)
\put(541,-543){\makebox(0,0)[lb]{\smash{{\SetFigFont{10}{12.0}{\rmdefault}{\mddefault}{\updefault}{\color[rgb]{0,0,0}$\beta_v$}%
}}}}
\end{picture}%

\caption {Given a crossing $v$ in the link diagram as on the left, we
construct a piece of the Heegaard surface $\Sigma$ on the right. This
piece contains four bits of alpha curves, shown in dashed lines, and
one beta curve $\beta_v.$}
\label {fig:inters}
\end {center}
\end {figure}

\begin {figure}
\begin {center}
\input {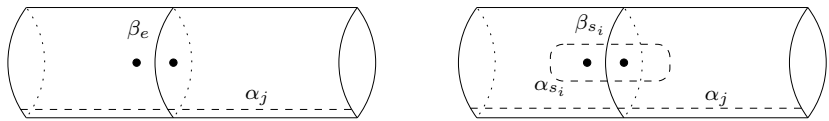}
\caption {A neighborhood of the distinguished edge $e$ (left) and a 
ladybug around some edge $s_i$ (right).}
\label {fig:lady}
\end {center}
\end {figure}

We denote by $\Sigma$ the boundary of a regular neighborhood of $D$ in 
$S^3,$ a surface of genus $g = c+1.$ To every region $A_r \ (r>0)$ we 
associate a curve $\alpha_r$ on $\Sigma,$ following the boundary of $A_r.$ 
To each crossing $v$ in $D$ we associate a curve $\beta_v$ on $\Sigma$ as 
indicated in Figure~\ref{fig:inters}. Furthermore, we introduce an extra 
curve $\beta_e$ which is the meridian of the knot, supported in a 
neighborhood of the distinguished edge $e.$ We also puncture the surface 
$\Sigma$ at two points on each side of $\beta_e,$ as shown on the 
left side of Figure~\ref{fig:lady}. 

Finally, for every edge $s_i, \ i=1, \dots, k,$ we introduce a ladybug, 
i.e. an additional pair of alpha-beta curves on $\Sigma$, as well as an 
additional pair of punctures. This type of configuration is shown on the 
right side of Figure~\ref{fig:lady}. The new beta curve $\beta_{s_i}$ is a 
meridian of the link, the two punctures lie on each of its two sides, and 
the new alpha curve $\alpha_{s_i}$ encircles the punctures.

In fact, we can think of the surface $\Sigma$ as being the union of 
several pieces, namely four-punctured spheres as in 
Figure~\ref{fig:inters} for each of the crossings, together with cylinders 
associated to $e$ and all the edges $s_i$ as in Figure~\ref{fig:lady}. 
Note that there could be several cylinders for the same edge.
The surface $\Sigma,$ together with the collections of curves and 
punctures, is an example of a multi-pointed Heegaard diagram for $S^3$ 
compatible with $L,$ in the sense of \cite[Definition 2.1]{MOS}. (However, 
unlike in the original sources, here we are only interested in the hat 
version of $HFK,$ and hence we do not need to distinguish between two 
different types of punctures.)

\begin {figure} \begin {center} \input {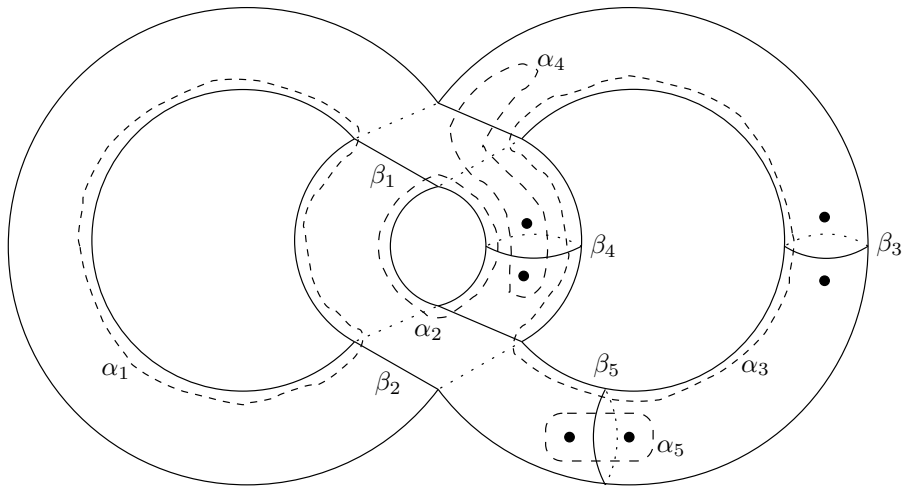} 
\caption {This is a 
special Heegaard diagram compatible with the Hopf link, with $g=3$ and 
$k=2.$ The beta curves $\beta_1$ and $\beta_2$ are associated to the two 
crossings, $\beta_3$ to the distinguished edge, while $\beta_4$ and 
$\beta_5$ are each part of a ladybug.  There are three alpha curves 
associated to planar bounded regions and two, $\alpha_4$ and $\alpha_5,$ 
which are parts of ladybugs. One tip of $\alpha_4$ is stretched to achieve 
admissibility.} 
\label {fig:hopf} \end {center} \end {figure}

For the purpose of defining Floer homology, we need to ensure that the 
Heegaard diagram is admissible in the sense of \cite[Definition 3.5]{OS5}. 
This condition can be achieved by isotoping the curves, cf. \cite{OS1}, 
\cite{OS5}. For example, one could stretch one tip of the alpha curve of 
each ladybug, and bring it close to the punctures associated to the 
distinguished edge $e.$ It is easy to see that the result is an admissible 
diagram; see Figure~\ref{fig:hopf} for an example. In general, we get a 
diagram with $k$ ladybugs, $g+k$ alpha curves, $g+k$ beta curves, and 
$2k+2$ punctures in this diagram. We denote by $\wsigma$ the complement of 
the punctures in the surface $\Sigma.$

We then consider the torus $\ta$ which is the product of all the alpha 
curves and the torus $\tb$ which is the product of all the beta
curves. We view $\ta$ and $\tb$ as totally real 
submanifolds of the symmetric product $\sym^{g+k}(\wsigma).$ Let 
$CF(\ta, 
\tb)$ be the vector space freely generated by the intersection points 
between $\ta$ and $\tb.$ One endows $CF(\ta,\tb)$ with the differential
$$ \del \xx = \sum_{\yy \in \ta \cap \tb} \sum_{\{\phi\in\pi_2({\mathbf
    x},{\mathbf y})\big|
\mu(\phi)=1 \}}
\#\left(\frac{{\mathcal
      M}(\phi)}{\mathbb R}\right) {\mathbf y}.$$
Here $\pi_2({\mathbf x},{\mathbf y})$ denotes the space of
homology classes of Whitney disks connecting ${\mathbf x}$ to
${\mathbf y}$ in $\wsigma$, ${\mathcal M}(\phi)$ denotes the moduli space 
of pseudo-holomorphic representatives of $\phi$ (with respect to a 
suitable almost complex structure), and $\mu(\phi)$ denotes its
formal dimension (Maslov index).

We can take the homology with respect to $\del,$ and obtain Floer homology 
groups $HF(\ta, \tb).$ According to \cite[Proposition 2.4]{MOS}, these 
are (up to a factor) the knot Floer homology groups of \cite{OS2} and 
\cite{R1}: 
$$ HF(\ta, \tb) \cong \hfk(L) \otimes V^{k-l}.$$

\section {Proof of Theorem~\ref{main}}
\label {sec:main}

Consider three links $L, L_0, L_1,$ with planar diagrams $D, D_0, D_1$ 
differing from each other as in the statement of Theorem~\ref{main}. Among 
them there is exactly one which has $m$ components, while the other two 
have only $m-1$ components. Without loss of generality, let us assume that 
$L$ has $m$ components. (The other two cases are similar.)

Pick a special Heegaard diagram for $L$ as in Section~\ref{sec:special}. 
We choose it to have the minimum possible number of ladybugs, namely 
$m-1.$ We denote the alpha and beta curves in the diagram by $\alpha_i$ 
and $\beta_i,$ respectively, with $i=1, \dots, n,$ where $n=g+m-1.$ We 
reserve the index $n$ for the beta curve $\beta=\beta_n$ associated to the 
particular crossing $v$ where $D$ differs from $D_0$ and $D_1.$

\begin {figure}
\begin {center}
\begin{picture}(0,0)%
\includegraphics{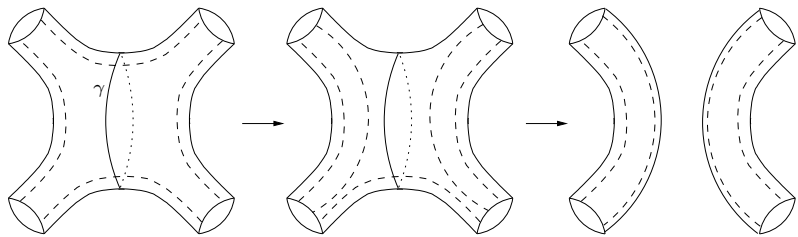}%
\end{picture}%
\setlength{\unitlength}{2881sp}%
\begingroup\makeatletter\ifx\SetFigFont\undefined%
\gdef\SetFigFont#1#2#3#4#5{%
  \reset@font\fontsize{#1}{#2pt}%
  \fontfamily{#3}\fontseries{#4}\fontshape{#5}%
  \selectfont}%
\fi\endgroup%
\begin{picture}(6793,1967)(837,-1716)
\put(1576,-488){\makebox(0,0)[lb]{\smash{{\SetFigFont{7}{8.4}{\rmdefault}{\mddefault}{\updefault}{\color[rgb]{0,0,0}$\gamma$}%
}}}}
\end{picture}%

\caption {A handleslide of the alpha curves, then a 
de-stabilization preceded by suitable handleslides of beta curves over 
$\gamma.$ The result is a special Heegaard diagram for the resolution $L_0.$}
\label {fig:moves}
\end {center}
\end {figure}

Let us call $\mathcal{P}$ the piece of the Heegaard diagram corresponding 
to the crossing $v$ as in Figure~\ref{fig:inters}. Topologically, 
$\mathcal{P}$ is a sphere with four disks removed. In $\mathcal{P},$ let 
us replace $\beta$ by the curve $\gamma$ pictured on the left side of 
Figure~\ref{fig:moves}. It can be checked directly that this gives a 
multi-pointed Heegaard diagram for the resolution $L_0.$ Alternatively, 
one can perform the moves shown in Figure~\ref{fig:moves} to arrive at a 
special Heegaard diagram for $L_0,$ of the type considered in the previous 
section.

\begin {figure}
\begin {center}
\input {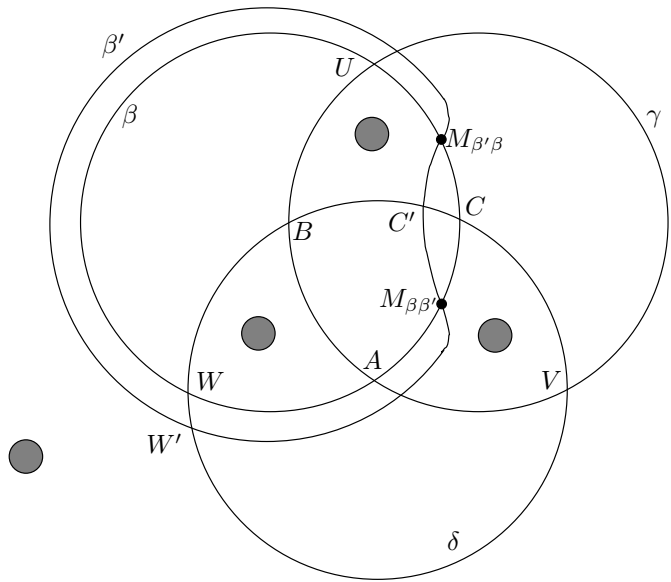}
\caption {We show various curves in the four-punctured sphere 
$\mathcal{P}$. The four gray disks correspond to tubes which 
link $\mathcal{P}$ to the rest of the 
Heegaard surface $\Sigma.$} \label {fig:triad}
\end {center}
\end {figure}

A similar construction can be used to obtain a Heegaard diagram for $L_1$: 
instead of the vertical curve $\gamma$ in $\mathcal{P}$, we need to 
consider a horizontal curve $\delta,$ separating the two upper boundaries 
of $\mathcal{P}$ from the two lower boundaries. If we identify $S^2$ with 
the plane together with a point at infinty, and $\mathcal{P}$ with the 
complement of four disks in $S^2$, then the curves $\beta, \gamma,$ and 
$\delta$ are arranged as in Figure~\ref{fig:triad}. Their intersection 
points are: $$ \beta \cap \gamma = \{A,U\} , \ \ \gamma \cap \delta= \{B, 
V \}, \ \ \delta \cap \beta = \{C, W \}.$$

Consider also a small perturbation $\beta'$ of the curve $\beta,$ such 
that the two intersect in two points. We denote by $C'$ and 
$W'$ the intersection points of $\beta'$ and $\delta$ which 
are close to $C$ and $W,$ respectively. Furthermore, for each $i=1, \dots, 
n-1,$ we choose small perturbations $\gamma_i, \delta_i, \beta_i'$ of the 
curve $\beta_i$, such that (for a fixed $i$) each intersects $\beta_i$ in 
two points, and any two of them intersect each other in two points as 
well. In general, if $\eta$ and $\eta'$ are (homologically nontrivial) 
isotopic curves in the punctured surface $\hat \Sigma$ intersecting in two 
points, then the Floer chain complex $CF(\eta, \eta')$ has rank two, being 
generated by the intersection points. We denote by $M_{\eta \eta'}$ the 
point which gives the top degree generator of $CF(\eta, \eta'),$ and by 
$M_{\eta'\eta}$ the other intersection point. For example, $\beta$ and 
$\beta'$ intersect in the points $M_{\beta' \beta}$ and $M_{\beta 
\beta'}$ shown in Figure~\ref{fig:triad}.

Set $ \ta = \alpha_1 \times \dots \times \alpha_n$ and
$$ \tb = \beta_1 \times \dots \times \beta_{n-1} \times \beta; \ \ \tg = 
\gamma_1 \times \dots \times \gamma_{n-1} \times \gamma; $$
$$ \td = \delta_1 \times \dots \times \delta_{n-1} \times \delta; \ \ \tbp 
= \beta'_1 \times \dots \times \beta'_{n-1} \times \beta'. $$

Viewing them as totally real tori in $\sym^n(\wsigma),$ we have
$$ HF(\ta, \tb) = HF(\ta, \tbp)= \hfk(L); $$
$$HF(\ta, \tg) = \hfk(L_0) \otimes V; \ HF(\ta, \td) = \hfk(L_1) \otimes 
V.$$

Thus, the exact triangle from the statement of Theorem~\ref{main} can be 
written as:
\begin {equation}
\label {abc}
\dots \To HF(\ta, \tb) \To HF(\ta, \tg) \To HF(\ta, \td) \To \dots.
\end {equation}

The strategy for proving (\ref{abc}) is the same as the one used by 
Ozsv\'ath and Szab\'o in \cite{OS4} for double-branched covers. Roughly, 
the principle behind the proof is the following: ``If the triangle counts 
are zero, and the quadrilateral counts are one (modulo two), then the 
exact triangle holds true.''

To make this precise, we need the following input from homological 
algebra, which was also used in \cite{OS3} and \cite{KMOS}:
\begin {lemma}
\label {hom}
Let $\{ (C_k, \del_k) \}_{k \in \zz/3\zz}$ be a collection of chain 
complexes over 
$\ff=\zz/2\zz$ and let
$ \{ f_k: C_k \to C_{k+1} \}_{k \in \zz/3\zz}$
be a collection of chain maps with the following properties:
\begin {enumerate}
\item { The composite $f_{k+1} \circ f_k: C_k \to C_{k+2}$ is 
chain-homotopic 
to zero, by a chain homotopy $H_k$:
$$ \del_{k+2} \circ H_k + H_k \circ \del_k = f_{k+1} \circ f_k;$$ }
\item {The sum
$$ \psi_k = f_{k+2} \circ H_k + H_{k+1} \circ f_k: C_k \to C_k$$
(which is a chain map) induces an isomorphism on homology.} 
\end {enumerate}
Then the sequence
$$ \begin {CD}
\cdots @>>> H_*(C_{k-1}) @>(f_{k-1})_*>> H_*(C_k) @>(f_k)_*>> 
H_*(C_{k+1}) @>>> \cdots 
\end {CD}
$$
is exact.
\end {lemma}

In our situation, we seek to apply Lemma~\ref{hom} to 
$$ C_0 = CF(\ta, \tb), \ \ C_1 = CF(\ta, \tg), \ \ C_2 = CF(\ta, \td).$$

The maps $f_k$ and $H_k$ will be particular examples of the following 
well-known natural maps in Floer homology, cf. \cite{deS}, \cite{FOOO}. 
For Lagrangians $T_0, T_1, \dots, T_s$ in a symplectic manifold $M$ (with 
good topological and analytic properties, e.g. transversality, no 
bubbling), the choice of a compatible almost complex structure produces 
maps 
$$ F_{T_0, \dots, T_s}: \bigotimes_{i=1}^s CF(T_{i-1}, T_i) \To 
CF(T_0, T_s)$$ 
$$ F_{T_0, \dots, T_s}(\mathbf{x}_1 \otimes \cdots \otimes 
\mathbf{x}_s) = \sum _{\mathbf{y}\in T_0 \cap T_s} 
\sum_{ \{\phi \in \pi_2(\mathbf{x}_1, \dots, \mathbf{x}_s, \mathbf{y}) | 
\mu(\phi)=0\} }
(\# \mathcal{M}(\phi)) 
\cdot 
\mathbf{y}.$$ 
Here $\pi_2(\mathbf{x}_1, \dots, \mathbf{x}_s, \mathbf{y})$ denotes the 
set of homotopy classes of Whitney $(s+1)$-gons in $M,$ with boundaries on 
$T_0, \dots, T_s$ and vertices $\mathbf{x}_1,\dots,\mathbf{x}_s, 
\mathbf{y};$ also, $\mu(\phi)$ is the Maslov index, and $
\mathcal{M}(\phi)$ is the moduli space of pseudo-holomorphic 
representatives of $\phi.$ (When $s=1,$ we ask for $\mu(\phi)=1$ instead 
of $0,$ and divide out the moduli space by the automorphism group $\rr$ 
before counting. The map $F_{T_0, T_1}$ is then the differential in the 
Floer complex.) These maps satisfy the generalized associativity 
relations: \begin {equation} \label {assoc} \sum_{0 \leq i < j \leq s} 
F_{T_0, \dots, T_{i-1}, T_i, T_j, \dots, T_s} \circ F_{T_i, \dots, T_j} 
=0. \end {equation}

The maps $F_{T_0, \dots, T_s}$ are also well-defined, and satisfy 
(\ref{assoc}), when the $T_i$'s are totally-real product tori in 
$\sym^n(\wsigma),$ cf. \cite{OS1}, \cite[Section 4.2]{OS4}.

Before applying these maps in our setting, let us introduce some more 
notation. Given one of the intersection points 
shown in Figure~\ref{fig:triad}, for example $A \in \beta \cap \gamma,$ we 
obtain a corresponding generator in the corresponding Floer chain complex 
$CF(\tb, \tg)$ by adjoining to $A$ the top degree intersection points 
$M_{\beta_i\gamma_i} \in \beta_i \cap \gamma_i.$ We denote the resulting 
generator by the respective lowercase letter in bold; for example: 
$$ \mathbf{a} = M_{\beta_1\gamma_1} \times M_{\beta_2\gamma_2}\times \dots 
\times M_{\beta_{n-1}\gamma_{n-1}}\times A \in 
CF(\tb, \tg);$$ 
$$\mathbf{u} = M_{\beta_1\gamma_1}\times M_{\beta_2\gamma_2}\times \dots 
\times M_{\beta_{n-1}\gamma_{n-1}} \times U \in
CF(\tb, \tg);$$
$$ \mathbf{b} = M_{\gamma_1\delta_1}\times M_{\gamma_2\delta_2} \times 
\dots \times M_{\gamma_{n-1}\delta_{n-1}} \times B 
\in CF(\tg, \td), \ \text{etc.} $$

We can now define the maps needed in Lemma~\ref{hom}. We choose
$$f_0: CF(\ta, \tb) \to CF(\ta, \tg), \ f_0(\mathbf{x}) = F_{\ta, \tb, 
\tg} \bigl(\mathbf{x} \otimes (\mathbf{a} + \mathbf{u} ) \bigr);$$ 
$$f_1: CF(\ta, \tg) \to CF(\ta, \td), \ f_1(\mathbf{x}) = F_{\ta, \tg,
\td} \bigl(\mathbf{x} \otimes (\mathbf{b} + \mathbf{v} ) \bigr);
$$
$$f_2: CF(\ta, \td) \to CF(\ta, \tb), \ f_2(\mathbf{x}) = F_{\ta, \td,
\tb} \bigl(\mathbf{x} \otimes (\mathbf{c} + \mathbf{w} ) \bigr);
$$
$$ H_0: CF(\ta, \tb) \to CF(\ta, \td), \ H_0(\mathbf{x}) = F_{\ta, \tb,
\tg, \td} \bigl(\mathbf{x} \otimes (\mathbf{a} + \mathbf{u} ) 
\otimes (\mathbf{b} + \mathbf{v} ) \bigr);$$
$$H_1: CF(\ta, \tg) \to CF(\ta, \tb), \ H_1(\mathbf{x}) = F_{\ta, \tg,
\td, \tb} \bigl(\mathbf{x} \otimes (\mathbf{b} + \mathbf{v} ) \otimes 
(\mathbf{c} + \mathbf{w} ) \bigr);$$
$$H_2: CF(\ta, \td) \to CF(\ta, \tg), \ H_2(\mathbf{x}) = F_{\ta, \td,
\tb, \tg} \bigl(\mathbf{x} \otimes (\mathbf{c} + \mathbf{w} )  
\otimes (\mathbf{a} + \mathbf{u}) \bigr).$$

\begin {lemma}
\label {one}
The maps $f_k \ (k\in \zz/3\zz)$ are chain maps.
\end {lemma}
\noindent {\bf Proof.} Equation (\ref{assoc}) for $s=2$ says that  if 
$\mathbf{a} + \mathbf{u}$ is a cycle in $CF(\tb, \tg),$ then the map $f_0$ 
commutes with the Floer differentials. We claim that both $\mathbf{a}$ and 
$\mathbf{u}$ are cycles. Note that if 
$u$ is a pseudo-holomorphic disk connecting $\mathbf{a}$ to some other 
intersection point $\mathbf{y} \in \tb \cap \tg,$ then $\mathbf{y}$ is an 
$n$-tuple of points, one of which is $y_n \in \{ A, U \}.$ Let 
$\mathcal{D}(u)$ be the domain in $\Sigma$ associated to 
$u,$ cf. \cite[Definition 2.13]{OS1}. 

Observe that if we delete from $\Sigma$ the $\beta$ and $\gamma$ curves, 
the only connected componets which do not contain punctures are the pairs 
of thin bigons joining $M_{\beta_i \gamma_i}$ to $M_{\gamma_i \beta_i},$ 
for each $i = 1, \dots, n-1.$ Since $u$ is a pseudo-holomorphic disk in 
$\sym^n(\wsigma),$ the domain $\mathcal{D}(u)$ cannot go over a puncture, 
hence it is a sum of those thin bigons (with some multiplicities). In 
particular, $y_n$ must be $A.$ 

If we impose the condition on $u$ to have Maslov index one, then 
$\mathcal{D}(u)$ is exactly one thin bigon. Each bigon has a holomorphic 
representative, hence contributes a term to the differential of 
$\mathbf{a}$. However, the bigons come in pairs, which means that $\del 
\mathbf{a} =0.$ Similarly, $\del \mathbf{u} =0.$ The cases of $f_1$ and 
$f_2$ are completely analogous. $\hfill \fin$ \medskip

\begin {lemma}
\label {two}
The maps $f_k$ and $H_k$ satisfy condition (1) in Lemma~\ref{hom}.
\end {lemma}
\noindent {\bf Proof.} We claim that
\begin {equation}
\label {abg1}
 F_{\tb, \tg, \td}(\mathbf{a} \otimes \mathbf{b})=F_{\tb, \tg, 
\td}(\mathbf{u} \otimes \mathbf{v})  = \mathbf{c}, \  
F_{\tb, \tg, \td}(\mathbf{a} \otimes 
\mathbf{v})=F_{\tb, \tg, \td}(\mathbf{u} \otimes \mathbf{b})  = 
\mathbf{w}. 
\end {equation} 
This follows from an inspection of the $\alpha \beta \gamma$ triangles in 
Figure~\ref{fig:triad}, together with the observation that the domains of 
pseudo-holomorphic triangles cannot go over the tubes, by an argument 
similar to the one in the proof of Lemma~\ref{one}.

Summing up the relations in (\ref{abg1}) we get:
\begin {equation}
\label {abg}
 F_{\tb, \tg, \td}\bigl((\mathbf{a}+\mathbf{u}) \otimes 
(\mathbf{b}+ \mathbf{v})\bigr)= 0.
\end {equation}

Using the associativity relation (\ref{assoc}) in the case $s=3,$ together 
with (\ref{abg}) and the fact that $\mathbf{a}+\mathbf{u}$ and 
$\mathbf{b}+\mathbf{v}$ are cycles, we obtain $\del_1 \circ H_0 + H_0 
\circ \del_0 = f_1 \circ f_0.$ The cases $k=1$ and $k=2$ are similar. 
$\hfill \fin$
\medskip

\begin {lemma}
\label {three}
The maps $f_k$ and $H_k$ satisfy condition (2) in Lemma~\ref{hom}.
\end {lemma}
\noindent {\bf Proof.} We will prove 
condition (2) for $k=0,$ as the other cases are completely analogous. 
Let 
$$\mathbf{m} = M_{\beta_1 \beta_1'} \times \dots \times  
M_{\beta_{n-1}\beta_{n-1}'} \times 
M_{\beta \beta'} \in CF(\tb, \tbp).$$
The map
$$ g: CF(\ta, \tb) \to CF(\ta, \tbp), \ g(\mathbf{x}) = F_{\ta, \tb, 
\tbp}(\mathbf{x} \otimes \mathbf{m})$$
is readily seen to be an isomorphism on homology. Observe also that the 
maps $g \circ f_2$ and $g \circ H_1$ are homotopy equivalent to $f_2'$ and 
$H_1',$ respectively, where 
$$f_2': CF(\ta, \td) \to CF(\ta, \tbp), \ f_2'(\mathbf{x}) = F_{\ta, \td,
\tbp} \bigl(\mathbf{x} \otimes (\mathbf{c'} + \mathbf{w'} ) \bigr);
$$
$$H_1': CF(\ta, \tg) \to CF(\ta, \tbp), \ H_1'(\mathbf{x}) = F_{\ta, \tg,
\td, \tbp} \bigl(\mathbf{x} \otimes (\mathbf{b} + \mathbf{v} ) \otimes
(\mathbf{c'} + \mathbf{w'} ) \bigr).$$

Therefore, the statement that $f_2 \circ H_0 + H_1 \circ f_0$ is an 
isomorphism on homology is equivalent to 
$$g_1 = f_2' \circ H_0 + H_1' \circ f_0 : CF(\ta, \tb) \to CF(\ta, \tbp)$$
being an isomorphism on homology.

Consider the map
$$ g_2: CF(\ta, \tb) \to CF(\ta, \tbp), \ g_2(\mathbf{x}) = F_{\ta, \tb,
\tbp} \bigl( \mathbf{x} \otimes \mathbf{\theta} \bigr),$$ 
where
$$ \mathbf{\theta} = F_{\tb, \tg, \td, \tbp} \bigl(  
(\mathbf{a}+\mathbf{u}) \otimes (\mathbf{b}+ \mathbf{v}) \otimes
(\mathbf{c'} + \mathbf{w'} )\bigr).$$

Finally, define $H: CF(\ta, \tb) \to CF (\ta, \tbp)$ by
$$ H(\mathbf{x}) = F_{\ta, \tb, \tg, \td, \tbp} \bigl (\mathbf{x}
\otimes (\mathbf{a}+\mathbf{u}) \otimes (\mathbf{b}+ \mathbf{v}) \otimes
(\mathbf{c'} + \mathbf{w'} )\bigr).$$

Let us apply the $s=4$ version of the associativity relation (\ref{assoc}) 
to the tori $\ta, \tb, \tg, \td, \tbp$ (in this order), evaluating 
all the summands at an element of the form $\mathbf{x}
\otimes (\mathbf{a}+\mathbf{u}) \otimes (\mathbf{b}+ \mathbf{v}) \otimes
(\mathbf{c'} + \mathbf{w'}).$ We get ten summands; three of them evaluate 
to zero because $\mathbf{a}+\mathbf{u}, \mathbf{b}+ \mathbf{v}$ and $
\mathbf{c'} + \mathbf{w'}$ are cycles, cf. Lemma~\ref{one}; two others 
evaluate to zero because
$$  F_{\tb, \tg, \td}\bigl((\mathbf{a}+\mathbf{u}) \otimes
(\mathbf{b}+ \mathbf{v})\bigr)=  
F_{\tg, \tb, \tbp}\bigl((\mathbf{b}+ \mathbf{v}) \otimes 
(\mathbf{c'} + \mathbf{w'})\bigr)= 0,$$
cf. Lemma~\ref{two}. The remaining five summands give the relation:
$$ g_1(\mathbf{x}) + g_2(\mathbf{x}) + (\del \circ H)(\mathbf{x}) + (H 
\circ \del)(\mathbf{x}) =0.$$

Therefore, $g_1$ and $g_2$ are chain homotopic. It suffices now to show 
that $g_2$ is an isomorphism on homology. At 
this time we need to make use of the homological (Maslov) grading on Floer 
complexes. We claim that
\begin {equation}
\label {ek}
\mathbf{\theta} = \mathbf{m} + \text{(lower degree terms)}.
\end {equation}
In other words, the claim is that the count of Maslov 
index zero, pseudo-holomorphic quadrilaterals having as vertices one of 
$\mathbf{a}$ and $\mathbf{u},$ one of $\mathbf{b}$ and $\mathbf{v},$ one 
of $\mathbf{c'}$ and $\mathbf{w'},$ as well as $\mathbf{m},$ is odd. This 
follows by inspecting the corresponding quadrilaterals in 
Figure~\ref{fig:triad}, and coupling them with quadrilaterals between 
$\beta_i, \gamma_i, \delta_i$ and $\beta'_i \ (i=1, \dots, n-1),$ as in 
the proof of Theorem 4.5 in \cite{OS4}. With the perturbation $\beta'$ of 
$\beta$ being chosen exactly as in Figure~\ref{fig:triad}, there is a 
unique useful (Maslov index zero) quadrilateral there, the one with 
vertices $A, B, C'$ and $M_{\beta \beta'}$ which is positioned right in the 
middle of the picture.

Equation (\ref{ek}) implies that $g_2 = g + $(lower degree terms), hence 
$g_2$ is an isomorphism on homology. $\hfill \fin$
\medskip

The proof of Theorem~\ref{main} is completed by putting together Lemmas 
\ref{hom}, \ref{one}, \ref{two} and \ref{three}.

\section {An example}

\begin {figure}
\begin {center}
\input {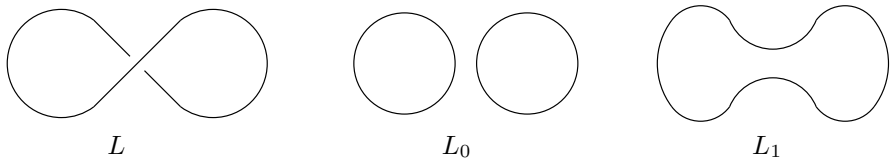}
\caption {A diagram of the unknot and its resolutions.}
\label {fig:unknot}
\end {center}
\end {figure}

\begin {figure}
\begin {center}
\input {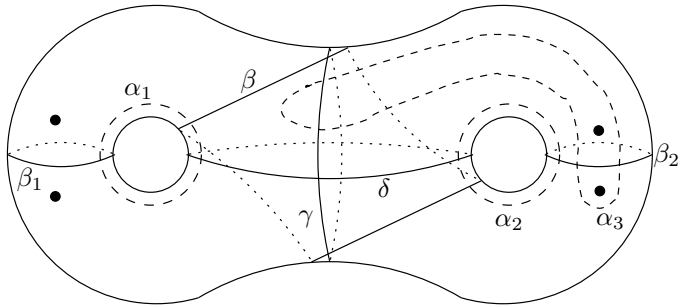}
\caption{The triples of curves $(\alpha_1, \alpha_2, \alpha_3)$ and 
$(\beta_1, \beta_2, \beta)$ give a Heegaard diagram for $L.$ To obtain 
Heegaard diagrams for $L_0$ and $L_1,$ replace the curve $\beta$ by 
$\gamma$ and $\delta,$ respectively.}
\label {fig:u1}
\end {center}
\end {figure}

We illustrate the unoriented skein exact triangle in the case of the 
links in Figure~\ref{fig:unknot}. The corresponding Heegaard 
diagrams are shown in Figure~\ref{fig:u1}. That picture can be simplified 
by handlesliding the curves $\beta$ and $\gamma$ over $\beta_1,$ 
handlesliding $\beta$ and $\delta$ over $\beta_2,$ and then de-stabilizing 
the pairs $(\alpha_1, \beta_1)$ and $(\alpha_2, \beta_2).$ 

\begin {figure}
\begin {center}
\input {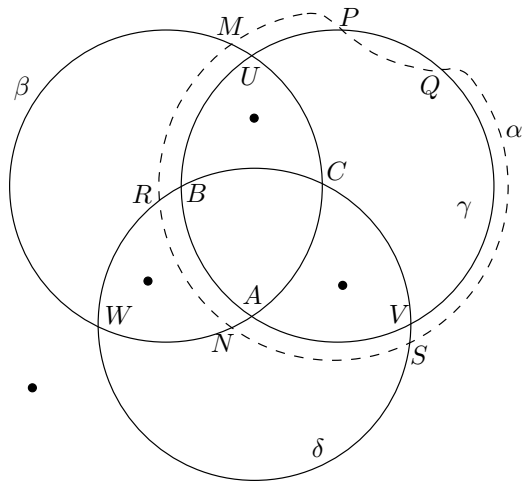}
\caption {Figure~\ref{fig:u1} after two de-stabilizations.}
\label {fig:u2}
\end {center}
\end {figure}

The result is shown in Figure~\ref{fig:u2}. The curve $\alpha = \alpha_3$
intersects each of the curves $\beta,$ $\gamma$ and $\delta$ in two 
points, so that each of the corresponding Floer homology groups are 
two-dimensional. The triangle (\ref{abc}) takes the form
$$ \begin {CD}
HF(\alpha, \beta) @>{f_0}>> HF(\alpha, \gamma) @>f_1>> HF(\alpha, 
\delta) @>f_2>> HF(\alpha, \beta).
\end {CD}
$$
The maps can be computed by inspecting the triangles in 
Figure~\ref{fig:u2}. For example, the $\alpha\delta\beta$ triangle with 
vertices $S, W$ and $N$ produces a summand of $N$ in $f_2(S).$ We obtain:
$$f_0(M) = Q; \ f_0(N) = Q; \ f_1(P) = R+S; \ f_1(Q) =  0;\ f_2(R) = M+N; 
\ f_2(S) = M+N.$$

\section {Quasi-alternating knots}
\label {sec:qa}

Let us recall the definition of the class $\mathcal{Q}$ of
quasi-alternating links from \cite{OS4}. The set $\mathcal{Q}$ is the 
smallest set of
links satisfying the following properties:
\begin {itemize}
\item {The unknot is in $\mathcal{Q}$;}
\item {If $L$ is a link which admits a projection with a crossing such
that
\begin {enumerate}
\item {Both resolutions $L_0, L_1 \in \mathcal{Q},$}
\item {$\det(L) = \det(L_0) + \det(L_1),$}
\end {enumerate}
then $L$ is in $\mathcal{Q}.$}
\end {itemize}

Note that if $L \in \mathcal{Q},$ then $\det(L) \geq 1,$ with equality if
and only if $L$ is the unknot. 
\medskip

\noindent \textbf{Proof of Corollary~\ref{cor}.} For any link $L,$ if we
use the bigradings on $\hfk(L)$ introduced in \cite{OS2}, \cite{R1}, then
the Euler characteristic of $\hfk(L)$ is the renormalized Alexander 
polynomial $(t^{-1/2} - t^{1/2})^l\cdot \Delta_L(t).$ The determinant of 
the link is $\det(L) =
|\Delta_L(-1)|,$ hence for every link we have:
$$ 2^{l-1} \cdot \det(L) \leq \rk \hfk(L).$$

We prove by induction on $\det(L)$ that a quasi-alternating link $L$ 
satisfies the inequality in the other direction, $ \rk \hfk(L) \leq 
2^{l-1} \cdot \det(L).$ Indeed, for the unknot we have $\rk \hfk=1,$ and 
then the inductive step follows readily from the definition of 
$\mathcal{Q}. \hfill \fin$

\begin {figure}
\begin {center}
\input {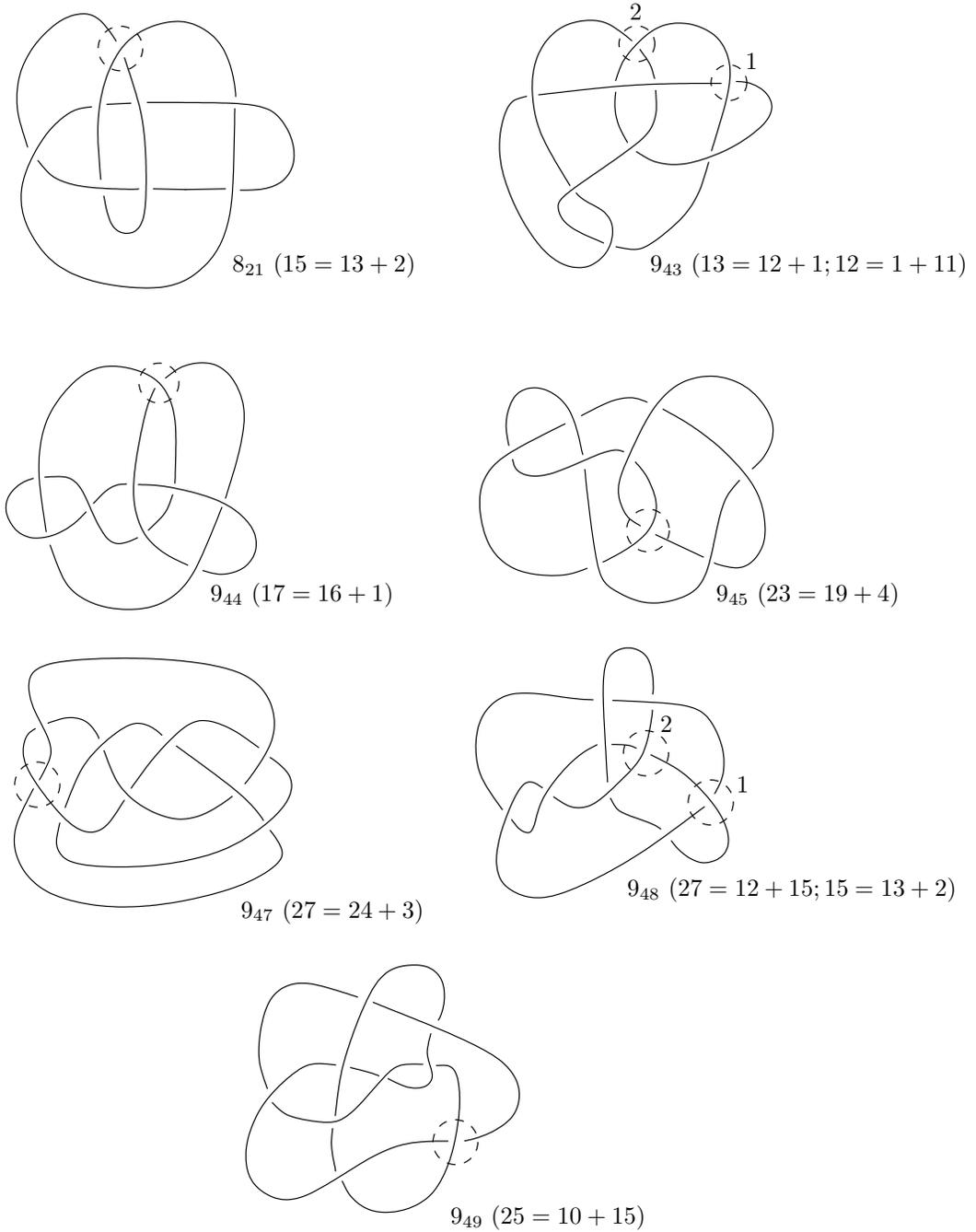}
\caption {Quasi-alternating knots with at most nine crossings that 
are not alternating. Five of them can be resolved at one crossing to 
produce two alternating links; we write the relation 
$\det(L)=\det(L_0)+\det(L_1)$ in parantheses. For the knots $9_{43}$ and 
$9_{48},$ one of the resolutions at the crossing numbered $1$ is 
alternating and the other is not; however, if we resolve the 
non-alternating one at the 
crossing numbered $2$, we obtain two alternating links.} \label 
{fig:quasi}
\end {center}
\end {figure}

\medskip We note that quasi-alternating knots are very frequent 
among small knots. Indeed, all but eleven of the prime knots with nine or 
less crossings are alternating, and alternating knots are 
quasi-alternating by \cite[Lemma 3.2]{OS4}. Furthermore, seven of the 
eleven non-alternating knots can be checked to be quasi-alternating; see 
Figure~\ref{fig:quasi}. In addition, John Baldwin informed the author 
that the knot $8_{20}$ is quasi-alternating as well. This leaves only the 
knots $8_{19}$ and $9_{42},$ 
which have $\rk \hfk (K) > \det (K) $ by \cite{OS2}, \cite{OS3.1}, and 
hence cannot be quasi-alternating, as well as 
$9_{46}$,  which the author does not know if it is 
quasi-alternating.\footnote{The published version of this paper contained  
a diagram exhibiting $9_{46}$ as a quasi-alternating knot. However, 
Abhijit Champanerkar informed the author that the diagram was incorrect.}

\end{document}